\documentclass[12pt,reqno]{article}
\usepackage{amsmath,amsfonts,amssymb,amsthm,indentfirst}

\catcode`\@=11

 \def\AMSTeXfeatures{\Plainheads 
   \let\current@vert=\AMS@vert}

 \def\Plainheads{\sh@ftdiam=0.05em
   \getlabeldims
   \let\vshaftfill=\plnvsolidfill
   \let\hshaftfill=\plnhsolidfill
   \let\th@rhead=\plnrhead
   \let\th@lhead=\plnlhead
   \let\th@dnhead=\plndnhead
   \let\th@uphead=\plnuphead}
 
 \def\glet{\global\let}

 \def\LaTeXfeatures{\catcode`\@=11
   \ifx\@clnwd\undefined \nol@g
      \input ltxcode.tex \dol@g \fi
   \ltxheads \let\current@vert=\new@vert
   \providelto \catcode`\@=\active}

 \def\nol@g{\def\wlog{\edef\garbage}}
 \def\dol@g{\let\wlog=\wl@g} \let\wl@g=\wlog
 \nol@g 

 \newbox\ltobox
 \def\providelto{{\setbox\z@=
   \hbox{$\to$}\minharrlen=\wd\z@
   \global\setbox\ltobox=\hbox{$\activeat>>>$}}
   \def\lto{\mathrel{\copy\ltobox}}}

 \def\ltxheads{\sh@ftdiam=\@wholewidth
   \getlabeldims
   \let\vshaftfill= \ltxvsolidfill
   \let\hshaftfill=\ltxhsolidfill
   \let\th@rhead=\ltxrhead
   \let\th@lhead=\ltxlhead
   \let\th@dnhead=\ltxdnhead
   \let\th@uphead=\ltxuphead}
 {\catcode`\@=\active
   \gdef@#1{\csname #1\string@at\endcsname}
   \glet\activeat=@}
 \def\def@#1{\expandafter\def\csname #1@at\endcsname}

 \def@>#1>#2>{\@rrow R{#1}{#2}}
 \def@<#1<#2<{\@rrow L{#1}{#2}}
 \def@ V#1V#2V{\@rrow V{#1}{#2}}
 \def@ A#1A#2A{\@rrow A{#1}{#2}}
 \def@/#1/#2/#3/{\@rrow{#1}{#2}{#3}}
 \def@.{\ifodd\row\ifmmode\noharrow
     \else\leavevmode.\spacefactor3000 \fi
   \else\novarrow\fi}
 \def@={\ifodd\row\harrow\hequalfill{}{}%
   \else\varrow\vequalfill{}{}\fi}
 \def@:#1{\ifx=#1\harrow\deffill{}{}%
   \else\leavevmode\null:#1\fi}
 \def@|{\current@vert}
  \def\AMS@vert{\varrow\vequalfill{}{}}
  \def\new@vert#1|#2|{\ifodd\row
   \let\nextarrow\vertexvarrow
   \else\let\nextarrow\varrow\fi
   \nextarrow\vshaftfill{#1}{#2}}
 \def@-{\ifmmode\let\next\hl@ne
   \else\let\next\AMSatdash \fi \next}
  \def\hl@ne#1-#2-{\harrow\hshaftfill{#1}{#2}}
  \def\AMSatdash{\let\next\relax\leavevmode
    \def\next@{\ifx\next-%
      \def\next-{\futurelet\next\nextii@}%
     \else\def\next{\hbox{-}}\fi\next}%
    \def\nextii@{\ifx\next-\def\next-{\hbox{---}}%
      \else\def\next{\hbox{--}}\fi\next}%
    \futurelet\next\next@}
 \def@(#1){\tweenarrows{#1}}
 \def@[#1]{\setsp@n#1\relax\activeat}
 \def\fiberbox{\hbox{$\vcenter{\hr@le\hbox{\vr@le
   \kern1ex\vbox{\kern1.2ex}\vr@le}\hr@le}$}}
  \def\hr@le{\hrule height \sh@ftdiam}
  \def\vr@le{\vrule width \sh@ftdiam}
 \def@+#1+#2+#3+{\ifodd\row \harrow{#1}{#2}{#3}%
   \else \varrow{#1}{#2}{#3}\fi}

 \def\Dnarrfill{\vequalfill\Dnhe@d}
 \def\Uparrfill{\Uphe@d\vequalfill}
 
 \def\ontofill{\rtarrfill\kern-0.3em 
   \th@rhead\kern 0.3em} 

 \def\rtarrfill{\hshaftfill\th@rhead}
 \def\ltarrfill{\th@lhead\hshaftfill}
 \def\dnarrfill{\vshaftfill\th@dnhead}
 \def\uparrfill{\th@uphead\vshaftfill}
 \def\hequalfill{\plnhfill=}
 \def\deffill{:\plnhfill=}
 \def\plnvextfill#1{\setbox\z@
   \hbox{\the\textfont3 #1}%
   \dimen@=\dp\z@\advance\dimen@\ht\z@
   \copy\z@ \kern-\dimen@ 
   \cleaders\copy\z@ \vfill
   \kern-\dimen@ 
   \box\z@}
 \def\plnhfill#1{$\m@th\mkern-1.5mu\mathord#1\mkern-6mu
    \cleaders\hbox{$\mkern-2mu\mathord#1\mkern-2mu$}\hfill
    \mkern-6mu\mathord#1\mkern-1.5mu$}
 \def\vequalfill{\plnvextfill{\char'167}}
 \def\plnvsolidfill{\plnvextfill{\char'077}}
 \def\plnhsolidfill{\plnhfill-}
 \def\ltxhsolidfill{\leaders\hrule height\topofshaft depth\botofshaft
   \hfill}
 \def\ltxvsolidfill{\leaders\vrule width\sh@ftdiam\vfill}
 \def\hdashfill{\hd@sh\wd@sh
   \xleaders \hbox{\wd@sh\hd@sh\wd@sh}\hfill
   \wd@sh\hd@sh}
 \def\vdashfill{\vd@sh\wd@sh
   \xleaders \vbox{\wd@sh\vd@sh\wd@sh}\vfill
   \wd@sh\vd@sh}
 \def\dashed{\ifinmeasureCD\else
    \ifodd\row\option{\let\hshaftfill=\hdashfill}%
   \else\option{\let\vshaftfill=\vdashfill}\fi\fi}

 \newdimen\CDstrutht  \newdimen\CDstrutdp
 \newdimen\CDstrutlen \CDstrutlen=\CDstrutht
   \advance\CDstrutlen by \CDstrutdp

 \def\CDstrut{\vrule
   height \ifnum\row=1 \z@\else\CDstrutht \fi
   depth \ifnum\row=\numrows \z@ \else\CDstrutdp \fi
   width\z@}

 \newdimen\CDarrsurr \CDarrsurr=0.375em
 \newdimen\CDdashlen
 \newdimen\CDvarrlen \CDvarrlen=1.5\baselineskip
 \newdimen\minharrlen 
  \setbox\z@\hbox{$\longrightarrow$} \minharrlen=\wd\z@
 \newdimen\minCDharrlen \minCDharrlen=2.5em 
\newdimen \minc@lwd
\def\findminc@lwd{\minc@lwd=2\CDarrsurr
  \advance\minc@lwd\minCDharrlen}

 \newdimen\sh@ftdiam

 \newdimen\labelsurr \labelsurr=1.25 em

\newcount\sp@ncnt \sp@ncnt=\@ne
\newcount\sp@ncnt@ \sp@ncnt@=\@ne
\newdimen\@rrwd \newdimen\@rrdp

 \def\adjustbot#1{\option{\advance\@rrdp#1\relax}}
 
\def\pushvertex#1{\global\p@shlen#1\relax
   \global\let\maybepush=\dopush}

 \newdimen\p@shlen \p@shlen=\z@

 \let\maybepush=\relax
 \def\dopush{\ifinmeasureCD 
   \advance\locdimen by -\p@shlen 
   \else\advance \@rrwd by -\p@shlen \fi 
   \global\let\maybepush=\relax \global\p@shlen=\z@\relax}

 \def\span@ne{\global\sp@ncnt=\@ne\relax}
 \def\setsp@n#1#2{\global\sp@ncnt=#1\relax
   \ifx\relax#2\relax\else\global\sp@ncnt@=#2\relax\fi}

 \def\plnrhead{\llap{$\rightarrow\mkern-1.5mu$}}
 \def\plnlhead{\rlap{$\mkern-1.5mu\leftarrow$}}

 \def\clap#1{\hbox to \z@{\hss #1\hss}}

 \def\plndnhead{\hbox{\the\textfont3 \char'171}}
 \def\plnuphead{\hbox{\the\textfont3 \char'170}}
 \def\Dnhe@d{\hbox{\the\textfont3 \char'177}}
 \def\Uphe@d{\hbox{\the\textfont3 \char'176}}

 \def\ltxrhead{\raise\@xisheight
   \llap{\smash{\@linefnt\@getrarrow(1,0)}}}
 \def\ltxlhead{\raise\@xisheight
   \rlap{\@linefnt\@getlarrow(-1,0)}}
 \def\ltxuphead{\setbox\z@=\rlap{%
   \kern\@halfwidth\@linefnt\char'66}%
   \copy\z@\kern-\ht\z@}
 \def\ltxdnhead{\setbox\z@=\rlap{%
   \kern\@halfwidth\@linefnt\char'77}%
   \ht\z@=\z@\box\z@}

 \def\wd@sh{\kern0.5\CDdashlen}
 \def\hd@sh{\vrule height\topofshaft depth\botofshaft
    width\CDdashlen}
 \def\vd@sh{\hrule height\CDdashlen
   depth\z@ width\sh@ftdiam}

\def\xylist{14{3434}13{2414}12{1723}%
  23{1413}34{1153}11{0867}43{0707}%
  32{0580}21{0414}31{0291}41{0}}
\newcount\tgtcnt@
\def\find@xyargs{\dimen@=\@rrdp
  \advance\dimen@ by \CDstrutlen
  \tgtcnt@=\dimen@ \dimen@=\@rrwd 
  \divide\dimen@ by \@m 
  \divide \tgtcnt@ by \dimen@ 
  \expandafter\testxy\xylist\relax
  \unitlength=\@xarg\@rrdp
  \divide\unitlength by\@yarg\relax}
\def\testxy#1#2#3{\ifnum\tgtcnt@>#3
    \@xarg=#1\relax \@yarg=#2\relax
    \let\next=\ignorerest
  \else\let\next\testxy\fi\next}
\def\ignorerest#1\relax{\relax}

\let\scalefactor=\@ne
\def\SWarrow{\find@xyargs\vector
  (-\@xarg,-\@yarg)\scalefactor\hskip-\wd\@linechar}
\def\NWarrow{\find@xyargs\vector
  (-\@xarg,\@yarg)\scalefactor\hskip-\wd\@linechar}
\def\NEarrow{\find@xyargs\vector
  (\@xarg,\@yarg)\scalefactor}
\def\SEarrow{\find@xyargs\vector
  (\@xarg,-\@yarg)\scalefactor}
\def\rightupline{\find@xyargs\@linelen=\scalefactor
     \unitlength\@sline}
\def\rightdownline{\find@xyargs\@yarg=-\@yarg\relax
     \@linelen=\scalefactor\unitlength\@sline}

\def\Sim{\ifodd\row\setbox\z@=\hbox{$\sim$}\dimen@=\ht\z@
 \advance\dimen@ by -\@xisheight
  \vbox{\box\z@\kern-\@xisheight\kern\dimen@}%
  \else\hbox{$\wr$}\fi}

\def\harrow#1#2#3{\inmeasureCDtrue\findminarrwd
  {#2}{#3}{\sp@ncnt\minharrlen}\inmeasureCDfalse\span@ne
  \mathrel{\hbox{\options\hplace{#1}\ulabel{#2}\dlabel{#3}}}}

\def\noharrow{\harrow\hfill{}{}}
\def\vertexvarrow#1#2#3{\findarrdp \@rrwd=\z@ \setsp@n\@ne\@ne
  \vbox to \z@{\kern-1.2\CDstrutht
  \rlap{\options\vplace{#1}\llabel{#2}\rlabel{#3}}\vss}}

\newif\ifinmeasureCD
\def\measurelabel#1{\setbox\z@
  \hbox{$\scriptstyle#1\kern\labelsurr$}%
  \ifdim\wd\z@>\@rrwd \@rrwd=\wd\z@\fi}
\def\findminarrwd#1#2#3{\@rrwd=#3\relax
   \measurelabel{#1}\measurelabel{#2}}
\def\findCDarrwd#1#2{\@rrwd=\minCDharrlen
   \measurelabel{#1}\measurelabel{#2}%
  }

\newcount\row \row=\@ne \newcount\col \col=\@ne 
 \newcount\numrows 
\numrows=\@ne
 \newcount\numcols
\newcount\arrspan \newdimen\vrtxhalfwd  \newbox\tempbox

\def\DANABUG{\advance\col by \@ne
 \@rrwd=\minCDharrlen
  \advance\@rrwd by \vrtxhalfwd
  \advance\@rrwd by \CDarrsurr
  \ifnum\col>\numcols \numcols=\col
     \newlocdimen{col\the\col}\locdimen=\@rrwd 
  \else \ifdim\@rrwd>\c@l \c@l=\@rrwd\fi\fi}

\def\drop#1\\{
  \findvrtxhalfsum\DANABUG\advance\row by 2 \measureinit}

\def\measureinit{\col=\@ne \vrtxhalfwd=-\CDarrsurr\arrspan=\@ne\@rrwd=\z@
   \setbox\tempbox=\hbox\bgroup$}
\def\measure{
  \let\harrow\measureCDarrow
  \let\CDCR=\measureCR 
   \findminc@lwd 
  \inmeasureCDtrue
  \row=\@ne \numcols=\z@ \measureinit}

\def\endmeasure{\findvrtxhalfsum\DANABUG
  \numrows=\row 
  \inmeasureCDfalse}

\def\newlocdimen#1{\advance\dimenc@unt by \@ne
  \ifnum\dimenc@unt<\insc@unt
     \else\errmessage{No room for the CD}\fi
  \dimendef\locdimen=\dimenc@unt
  \expandafter\dimendef\csname#1\endcsname=\dimenc@unt}

 \def\r@wc@l{\csname row\the\row col\the\col\endcsname}
 \def\c@l{\csname col\the\col\endcsname}

 \def\findvrtxhalfsum{$\egroup
  \newlocdimen{row\the\row col\the\col}
  \locdimen=\vrtxhalfwd 
  \vrtxhalfwd=0.5\wd\tempbox 
  \advance\vrtxhalfwd by \CDarrsurr
  \advance\locdimen by \vrtxhalfwd 
  \advance\@rrwd by \locdimen 
  \maybepush
  \divide\@rrwd by \arrspan\relax
  \ifdim\@rrwd<\minc@lwd
    \ifnum\col>\@ne \@rrwd=\minc@lwd\fi \fi
  \loop 
    \ifnum\col>\numcols \numcols=\col
       \newlocdimen{col\the\col}
       \locdimen=\@rrwd 
    \else \ifdim\@rrwd>\c@l \c@l=\@rrwd\fi \fi
   \ifnum\arrspan>\@ne
      \advance\arrspan by -1 \advance\col by \@ne
  \repeat }

 \def\measureCDarrow#1#2#3{\findvrtxhalfsum
   \arrspan=\sp@ncnt\relax\global\sp@ncnt=1\relax
   \advance\col by \@ne
   \findCDarrwd{#2}{#3}%
   \setbox\tempbox=\hbox\bgroup$}

 \newcount\dr@tn \dr@tn=\z@
 \def\locate#1:#2{\ifinmeasureCD\else
   \count@=-#1
   \multiply\count@ by 2
   \advance\count@ by #2
   \dimen@=\count@\@rrwd
   \ifnum\dr@tn=\@ne\relax \else\dimen@=-\dimen@ \fi
   \dimen@i=\@rrdp
   \ifnum\dr@tn>\z@\advance\dimen@i by \CDstrutlen \fi
   \dimen@i=\count@\dimen@i
   \count@=#2 \multiply\count@ by 2
   \divide\dimen@ by \count@
   \divide\dimen@i by \count@
   \lift\dimen@i\nudge\dimen@\fi}

\def\betweenCDrows{\advance\row by \@ne \col=\@ne
\options}

\def\hbegin{\hbox\bgroup\kern\c@l \kern-\r@wc@l$}
\def\hend{$\glet\maybepush\relax \CDstrut\egroup}
\def\vbegin{\setbox\tempbox=\hbox\bgroup$}
\def\vend{$\egroup\ht\tempbox=\z@\dp\tempbox\CDvarrlen
  \box\tempbox}
\def\setCD{\let\harrow=\setCDarrow
  \let\CDCR=\setCR 
  \row=\@ne \col=\@ne \hbegin}
\let\endsetCD=\hend 

\def\findarrwd{\@rrwd=\z@ \count@=\col \advance\count@ by\sp@ncnt
  \loop\ifnum\count@>\col \advance\count@ by -1
      \advance\@rrwd by\csname col\the\count@\endcsname\repeat}
\def\setCDarrow#1#2#3{\kern\CDarrsurr\advance\col by \@ne
  \findarrwd \advance\@rrwd by -\r@wc@l  
  \@rrdp=\z@ 
  \maybepush
  \advance\col by -\@ne \advance\col by \sp@ncnt \span@ne
  \hbox to \@rrwd{\options
   \@rrwd=\scalefactor\@rrwd\hss
   \hplace{#1}\ulabel{#2}\dlabel{#3}\hss}%
   \kern\CDarrsurr}

\newdimen\labspacei 
\newdimen\labspaceii 

\newdimen\@xisheight
  \@xisheight=\the\fontdimen22\textfont2
\newdimen\labelskip
  \labelskip=\the\fontdimen10\textfont3 
\newdimen\topofshaft
\newdimen\botofshaft
\newdimen\botofulabel
\newdimen\topofdlabel
\def\getlabeldims{
  \topofshaft=0.5\sh@ftdiam
  \botofshaft=\topofshaft
  \advance\topofshaft by \@xisheight  
  \advance\botofshaft by -\@xisheight  
  \botofulabel=\topofshaft
  \advance\botofulabel by \labelskip
  \topofdlabel=\botofshaft
  \advance\topofdlabel by \labelskip}

\def\ulabel{\ifnum\row=\@ne\let\next\ulabeli
   \else\let\next\ulabellap\fi\next}
\def\ulabeli#1{\vbox{
  \clap{\kern-\@rrwd$\scriptstyle#1$}%
  \kern\botofulabel}\maybeoffset}
\def\ulabellap#1{\vbox to \z@{\vss
  \clap{\kern-\@rrwd$\scriptstyle#1$}%
  \kern\botofulabel}\maybeoffset}
\def\dlabel{\ifnum\row=\numrows\let\next\dlabeli
   \else\let\next\dlabellap\fi\next}
\def\dlabeli#1{\vtop{\kern\topofdlabel
  \clap{\kern-\@rrwd$\scriptstyle#1$}%
  }\maybeoffset}
\def\dlabellap#1{\vbox to \z@{\kern\topofdlabel
  \clap{\kern-\@rrwd$\scriptstyle#1$}%
  \vss}\maybeoffset}
\def\rlabel#1{\vbox to \z@{\vss
  \rlap{\kern\labelskip$\scriptstyle#1$}%
  \vss\kern-\@rrdp}\maybeoffset}
\def\llabel#1{\vbox to \z@{\vss
  \llap{$\scriptstyle#1$\kern\labelskip}%
  \vss\kern-\@rrdp}\maybeoffset}
\def\swlabel#1{\vtop{\kern0.5\@rrdp
  \llap{$\scriptstyle#1$\kern\labelskip\kern-0.5\@rrwd}
  }\maybeoffset}
\def\nwlabel#1{\vbox{
  \llap{$\scriptstyle#1$\kern\labelskip\kern-0.5\@rrwd}%
  \kern-0.5\@rrdp}\maybeoffset}
\def\selabel#1{\vtop{\kern0.5\@rrdp
  \rlap{\kern0.5\@rrwd\kern\labelskip$\scriptstyle#1$}%
  }\maybeoffset}
\def\nelabel#1{\vbox{
  \rlap{\kern0.5\@rrwd\kern\labelskip$\scriptstyle#1$}%
  \kern-0.5\@rrdp}\maybeoffset}
\def\cplace#1{\vbox to \z@{\vss
  \clap{$#1$\kern-\@rrwd}%
  \kern-\@rrdp\vss}\maybeoffset}
\def\hplace#1{\hbox to \@rrwd{#1}\maybeoffset}
\def\vplace#1{\clap{\vbox to \z@{#1\kern-\@rrdp}}\maybeoffset}

\newdimen\nudgeamount \nudgeamount=\z@
\newdimen\liftamount \liftamount=\z@
\let\maybeoffset\relax
\newbox\offsetbox \newdimen\lastheight
\def\dooffset{
  \setbox\offsetbox=\lastbox \lastheight=\ht\offsetbox 
  \setbox\offsetbox=\vbox{\kern-\liftamount\box\offsetbox}%
  \ht\offsetbox=\lastheight
  \kern\nudgeamount\box\offsetbox\kern-\nudgeamount
  \global\nudgeamount=\z@ \global\liftamount=\z@
  \glet\maybeoffset=\relax}
\def\nudge#1{\ifinmeasureCD\else
  \global\advance\nudgeamount#1\relax
  \global\let\maybeoffset\dooffset\fi}
\def\lift#1{\ifinmeasureCD\else
  \global\advance\liftamount#1\relax
  \global\let\maybeoffset\dooffset\fi}

\def\findarrdp{\@rrdp=\CDvarrlen
  \ifnum\sp@ncnt@>1
    \advance\@rrdp by \CDstrutlen
    \multiply\@rrdp by \sp@ncnt@
    \advance\@rrdp by -\CDstrutlen \fi
 }

\def\varrow#1#2#3{\ifnum\sp@ncnt>\@ne 
     \sp@ncnt@=\sp@ncnt\relax\fi
  \findarrdp \@rrwd=\z@ 
  \kern\c@l
   \hbox to \z@{\options
   \@rrdp=\scalefactor\@rrdp
    \hss\vplace{#1}\llabel{#2}\rlabel{#3}\hss}%
  \global\advance\col by \@ne \setsp@n\@ne\@ne
  }

\def\novarrow{\varrow\vfill{}{}}

\def\tweenarrows#1{\findarrwd \findarrdp \setsp@n\@ne\@ne
  \rlap{\options\cplace{#1}}}

\def\usarrow #1#2#3{\dr@tn=\@ne
  \findarrwd \findarrdp \setsp@n\@ne\@ne 
  \rlap{\options\cplace{#1}\nwlabel{#2}\selabel{#3}}%
  \dr@tn=\z@}
\def\dsarrow #1#2#3{\dr@tn=\tw@
  \findarrwd \findarrdp \setsp@n\@ne\@ne 
  \rlap{\options\cplace{#1}\swlabel{#2}\nelabel{#3}}%
  \dr@tn=\z@}
 \def\@rrow#1{\csname #1@rrow\endcsname}
 \def\R@rrow{\harrow \rtarrfill}
 \def\L@rrow{\harrow \ltarrfill}
 \def\V@rrow{\varrow \dnarrfill}
 \def\A@rrow{\varrow \uparrfill}
 \def\SE@rrow{\dsarrow \SEarrow}
 \def\NW@rrow{\dsarrow \NWarrow}
 \def\SW@rrow{\usarrow \SWarrow}
 \def\NE@rrow{\usarrow \NEarrow}
 \def\DS@rrow{\dsarrow \dnslope}
 \def\US@rrow{\usarrow \upslope}
 \def\upslope{\find@xyargs
       \@linelen=\unitlength\@sline}
 \def\dnslope{\find@xyargs\@yarg=-\@yarg\relax
       \@linelen=\unitlength\@sline}

\newtoks\optionlist 
\optionlist={}
\let\options\relax
\def\dooptions{\the\optionlist\global\optionlist={}%
  \glet\options=\relax}
\def\option#1{\ifinmeasureCD\else
  \glet\options=\dooptions
  \global\optionlist=\expandafter{\the\optionlist\relax#1}\fi}
\def\wider#1{\ifinmeasureCD\else
   \option{\advance\@rrwd by #1}\fi}
\def\deeper#1{\ifinmeasureCD\else
   \option{\advance\@rrdp by #1}\fi}

{\def\\{\global\let\sptoken= }\\ }

\def\CR{\futurelet\nexttok\testCR}
\def\testCR{\ifx\nexttok\sptoken
   \let\next\eatspaceCR\else\let\next\CDCR\fi\next}
\def\eatspaceCR#1 {\CR}
\def\measureCR{\ifx\nexttok\endmeasure\let\nextCR\relax
    \else\let\nextCR\drop\fi\nextCR}
\def\setCR{\ifodd\row
  \ifx\nexttok\endsetCD\else\hend\betweenCDrows\vbegin\fi
  \else\vend\betweenCDrows\hbegin\fi}

\countdef\dimenc@unt=11
\def\CD#1\endCD{
   \begingroup\let\\=\CR
  \m@th\offinterlineskip
   \measure#1\endmeasure\null\,\vcenter{\setCD#1\endsetCD}\,
   \endgroup
    }

\ifx\@clnwd\undefined \nol@g\else\catcode`\ =14\relax\fi
 \font\@linefnt=line10 
 \newcount\@tempcnta
 \newcount\@tempcntb
 \newdimen\@tempdima
 \newdimen\@tempdimb
 \newdimen\@wholewidth
 \newdimen\@halfwidth
   \@wholewidth\fontdimen8\@linefnt \@halfwidth .5\@wholewidth
 \newdimen\unitlength
 \newcount\@xarg
 \newcount\@yarg
 \newcount\@yyarg
 \newbox\@linechar
 \newdimen\@linelen
 \newdimen\@clnwd
 \newdimen\@clnht
 \newif\if@negarg
 
 \def\@whilenoop#1{}

 \def\@whiledim#1\do #2{\ifdim #1\relax#2\@iwhiledim{#1\relax#2}\fi}

 \def\@iwhiledim#1{\ifdim #1\let\@nextwhile=\@iwhiledim 
         \else\let\@nextwhile=\@whilenoop\fi\@nextwhile{#1}}

 \def\@sline{\ifnum\@xarg< 0 \@negargtrue \@xarg -\@xarg \@yyarg -\@yarg
   \else \@negargfalse \@yyarg \@yarg \fi
 \ifnum \@yyarg >0 \@tempcnta\@yyarg \else \@tempcnta -\@yyarg \fi
 \ifnum\@tempcnta>6 \@badlinearg\@tempcnta0 \fi
 \ifnum\@xarg>6 \@badlinearg\@xarg 1 \fi
 \setbox\@linechar\hbox{\@linefnt\@getlinechar(\@xarg,\@yyarg)}%
 \ifnum \@yarg >0 \let\@upordown\raise \@clnht\z@
    \else\let\@upordown\lower \@clnht \ht\@linechar\fi
 \@clnwd=\wd\@linechar
 \if@negarg \hskip -\wd\@linechar \def\@tempa{\hskip -2\wd\@linechar}\else
      \let\@tempa\relax \fi
 \@whiledim \@clnwd <\@linelen \do
   {\@upordown\@clnht\copy\@linechar
    \@tempa
    \advance\@clnht \ht\@linechar
    \advance\@clnwd \wd\@linechar}%
 \advance\@clnht -\ht\@linechar
 \advance\@clnwd -\wd\@linechar
 \@tempdima\@linelen\advance\@tempdima -\@clnwd
 \@tempdimb\@tempdima\advance\@tempdimb -\wd\@linechar
 \if@negarg \hskip -\@tempdimb \else \hskip \@tempdimb \fi
 \multiply\@tempdima \@m
 \@tempcnta \@tempdima \@tempdima \wd\@linechar \divide\@tempcnta \@tempdima
 \@tempdima \ht\@linechar \multiply\@tempdima \@tempcnta
 \divide\@tempdima \@m
 \advance\@clnht \@tempdima
 \ifdim \@linelen <\wd\@linechar
    \hskip \wd\@linechar
   \else\@upordown\@clnht\copy\@linechar\fi}
 
 \def\@getlinechar(#1,#2){\@tempcnta#1\relax\multiply\@tempcnta 8
 \advance\@tempcnta -9 \ifnum #2>0 \advance\@tempcnta #2\relax\else
 \advance\@tempcnta -#2\relax\advance\@tempcnta 64 \fi
 \char\@tempcnta}
 
 \def\vector(#1,#2)#3{\@xarg #1\relax \@yarg #2\relax
 \@tempcnta \ifnum\@xarg<0 -\@xarg\else\@xarg\fi
 \ifnum\@tempcnta<5\relax
 \@linelen=#3\unitlength
 \ifnum\@xarg =0 \@vvector 
   \else \ifnum\@yarg =0 \@hvector \else \@svector\fi
 \fi
 \else\@badlinearg\fi}
 
 \def\@svector{\@sline
 \@tempcnta\@yarg \ifnum\@tempcnta <0 \@tempcnta=-\@tempcnta\fi
 \ifnum\@tempcnta <5
   \hskip -\wd\@linechar
   \@upordown\@clnht \hbox{\@linefnt  \if@negarg 
   \@getlarrow(\@xarg,\@yyarg) \else \@getrarrow(\@xarg,\@yyarg) \fi}%
 \else\@badlinearg\fi}
 
 \def\@getlarrow(#1,#2){\ifnum #2 =\z@ \@tempcnta='33\else
 \@tempcnta=#1\relax\multiply\@tempcnta \sixt@@n \advance\@tempcnta
 -9 \@tempcntb=#2\relax\multiply\@tempcntb \tw@
 \ifnum \@tempcntb >0 \advance\@tempcnta \@tempcntb\relax
 \else\advance\@tempcnta -\@tempcntb\advance\@tempcnta 64
 \fi\fi\char\@tempcnta}
 
 \def\@getrarrow(#1,#2){\@tempcntb=#2\relax
 \ifnum\@tempcntb < 0 \@tempcntb=-\@tempcntb\relax\fi
 \ifcase \@tempcntb\relax \@tempcnta='55 \or 
 \ifnum #1<3 \@tempcnta=#1\relax\multiply\@tempcnta
 24 \advance\@tempcnta -6 \else \ifnum #1=3 \@tempcnta=49
 \else\@tempcnta=58 \fi\fi\or 
 \ifnum #1<3 \@tempcnta=#1\relax\multiply\@tempcnta
 24 \advance\@tempcnta -3 \else \@tempcnta=51\fi\or 
 \@tempcnta=#1\relax\multiply\@tempcnta
 \sixt@@n \advance\@tempcnta -\tw@ \else
 \@tempcnta=#1\relax\multiply\@tempcnta
 \sixt@@n \advance\@tempcnta 7 \fi\ifnum #2<0 \advance\@tempcnta 64 \fi
 \char\@tempcnta}
\catcode`\ =10

\dol@g 
\catcode`\@=\active
\LaTeXfeatures

\newtheorem{thm}{Theorem} \newtheorem{prob}{Problem}

\theoremstyle{definition} 
\newtheorem*{rem*}{Remark}

           \def\int{\textrm{int}}

\numberwithin{equation}{section}
\date{}

\begin{document}

\title{Action-type axiomatizable classes of group representations}
\author{Aleko Gvaramia\qquad Boris Plotkin\\
 Abhazia State University \  Hebrew University\\
 \ \ \  Suhum, Abhazia\qquad  \ Jerusalem, Israel}

\maketitle
\begin{abstract}

The paper adjoins the book \cite{PV} and turns to be, in a sense,
its continuation. In the book the varieties of representations had
been considered. In the matter of fact, the varieties under
consideration are action-type varieties. This paper studies other
classes of representations,  axiomatizable in the special
action-type logic.

\end{abstract}

\section{Background (preliminaries)}

\subsection{Variety $Rep-K$}

Given a commutative with the unit ring $K$, consider a variety
which is also a category, denoted by $Rep-K$. The algebras of this
variety are two-sorted algebras.  These algebras are pairs, or
representations, $(V,G)$, where $V$ is $K$-module and $G$ a group,
acting in $V$. The action is treated as an operation in the
two-sorted algebra $(V,G)$.

The axioms are as follows:

\begin{itemize}
\item{ 1. the mapping $a \to a \circ  g$ is a $K$-linear mapping in
$V$,}

\item{ 2. $(a \circ g_1) \circ g_2 = a \circ g_1 g_2$,}

\item{ 3. $a \circ 1 = a$.}
\end{itemize}

Here $1$ is the unit in $G$, $a \in V, \ \ g \in G$. The operation
$\circ$ recovers the representation $\rho: G \to Aut V$, which is
sometimes identified with the algebra $(V, G)$, i.e., $\rho = (V,
G)$.

The pointed identities together with the identities of groups and
$K$-modules determine the variety $Rep-K$. As usual, this variety
is also a category. Its morphisms are two-sorted homomorphisms

$$
\mu = (\alpha, \beta): (V,G) \to (V', G'),
$$
where $\alpha: V \to V'$ a $K$-homomorphism of modules, $\beta: G
\to G'$ is a homomorphism of groups, and $(a\circ g)^\alpha =
a^\alpha \circ g^\beta$. We have: $Ker \mu = (Ker \alpha, Ker
\beta)$. Take $Ker \alpha = V_0, \ \ Ker \beta = H$. Then the pair
$(V_0, H)$ is a congruence of the representation $(V,G)$ in the
following sense: $V_0$ is a submodule in $V$, invariant under the
action of the group $G$, $H$ is a normal subgroup in $G$, acting
trivially in $V / V_0$. We have $(V,G) / Ker \mu = (V / V_0, G /
H)$. This gives the theorem on homomorphisms.

Let us consider a free representation $W = W(X,Y)$ Here the pair
$(X,Y)$ is a two-sorted set. We have: $W(X,Y) = (XKF(Y), F(Y))$,
where $ F(Y)=F$ is a free group over the set $Y$, $XKF = \Phi$ is
a free $KF$-module over the set $X$ and $KF$ is a group algebra.

The elements of $\Phi$ have the form $w=x_1 u_1 + ... + x_n u_n$,
$u_i \in KF$; the elements of $F$ are written as $f=f(y_1, ... ,
y_m)$. The action $\circ$ is defined by the rule:
$$
w \circ f = w f = x_1(u_1 f) + ... + x_n(u_n f).
$$

It is easy to understand that the mappings $\alpha: X \to V$ and
$\beta : Y \to G$ determine the homomorphism $\mu = (\alpha,
\beta): (\Phi, F) = W \to (V,G)$.

Consider further the sets $Hom(W, \rho) = Hom(W, (V,G))$. In the
situation of the finite $X=\{ x_1, ... , x_n\}$ and $Y=\{ y_1, ...
, y_m\}$ the sets $Hom(W, \rho)$ are treated as affine spaces.
There is a natural bijection $Hom(W,\rho) \to V^{(n)} \times
G^{(m)}$. A point $$((\alpha (x_1), ... , (\alpha (x_n)),
(\beta(y_1), ... , \beta(y_m))) = ((a_1, ... ,a_n), (g_1, ... ,
g_m)),$$ $a_i \in V, \ g_k \in G$, corresponds to the homomorphism
$\mu = (\alpha, \beta): (\Phi, F) \to (V,G)$.

\subsection{Logic in $Rep-K$.}

Given $X$ and $Y$, consider a signature $$L = L_{X,Y} = \{ \vee,
\wedge, \neg, \exists x, \exists y, x \in X, y \in Y \}.$$ For the
algebra $W=W(X,Y)$ we consider equalities $w\equiv 0$ and $f\equiv
1$. We treat these equalities as logical formulas. The equalities
of the first kind we call {\it action-type equalities}, while the
equalities of the second type we call {\it group equalities}.
Denote by $LW = L_{X,Y} W(X,Y)$ the absolutely free $L$-algebra
over the set of all $W$-equalities. Its elements we call formulas
(elementary, first order formulas) over the free representation
$W=(\Phi, F)$.

Let us consider an example of $L$-algebra. Take the set
$Hom(W,(V,G))$. Let $Set(W,(V,G))$ be a system of all subsets of
$Hom(W,(V,G))$. It is clear, that boolean operations are defined
in $Set(W,(V,G))$ and it is a boolean algebra. Let us define also
quantifiers. Let $A$ be a subset of $Hom(W,(V,G))$. We set:
$\mu(\alpha, \beta) \in \exists x A$, if there exists $\nu =
(\alpha ', \beta)$ such that $\nu \in A$ and $\alpha (x_1) =
\alpha '(x_1)$ for $x_1 \neq x$. Analogously we define $\exists y
A$. Here $\exists x $ and $\exists y$ are quantifiers of the
boolean algebra $Set(W,(V,G))$ in the sense of the following
definition.

The quantifier $\exists$ of a boolean algebra $B$ is a mapping
$\exists : B \to B$ with the properties:

\begin{itemize}

\item{ 1. $\exists 0 = 0$,}

\item{ 2. $ a < \exists a$,}

\item{ 3. $\exists(a \wedge \exists b) = \exists a \wedge \exists
b$.}

\end{itemize}

Here 0 is a zero in $B$, and $a,b \in B$. Define now a canonical
homomorphism of $L$-algebras $$Val = Val^W _{(V,G)} : LW \to
Set(W,(V,G)).$$ We set: $\mu= (\alpha, \beta) \in Val (w \equiv
0)$, if $w^\alpha = 0$ in $V$, $\mu= (\alpha, \beta) \in Val (f
\equiv 1)$, if $f^\beta = 1$ in $G$. Since $LW$ is free, this
determines $Val(u)$ for every formula $u \in LW$. Here the set of
homomorphisms $Val(u)$ is called {\it the value of the formula $u$
in the representation}  $(V,G)$. If $Val(u) = Hom(W,(V,G))$, then
it means that the formula $u$ holds in the representation $(V,G)$.

Define further separately the logic of action (action-type logic).
It is generated by equalities of action $w\equiv 0$ in the
signature $L_X=\{ \vee, \wedge, \neg, \exists x, x \in X \}$.
There are no quantifiers $\exists y$ and equalities $f \equiv 1$.
Denote this logic by $L_X W$.

\subsection{Classes of representations}

Consider classes  $\mathfrak{X}$ of representations  $(V,G)$ in
$Rep-K$ and, simultaneously, the sets of formulas $T$ in the logic
$LW$ with countable $X$ and $Y$. We establish Galois
correspondence
$$
T = \mathfrak X^* = \{ u \in LW | \ u \ \textrm{holds in every}\
(V,G)\in \mathfrak X \}
$$
$$
\mathfrak X = T^* = \{ (V,G) | \quad \textrm{every formula from }\
T \ \textrm{holds in} \ (V,G)\}.
$$
 A class $\frak X$ of the kind $\mathfrak X=T^*$ is called {\it axiomatizable one}.
Here $T$ is an arbitrary set. If $T$ consists of action-type
formulas, then $\mathfrak X=T^*$ is an  {\it action-type
axiomatizable class}.

Let us distinguish special cases.

1. $T$ consists of equalities. Then $\mathfrak X=T^*$ is a variety
of representations. If $T$ consists of action-type equalities,
then $\mathfrak X=T^*$ is {\it an action-type variety}. Such
$\mathfrak X$ is given also by the formulas of the form $x \circ u
\equiv 0$, $u \in KF$.

2. Formulas of the kind $u_1 \vee u_2 \vee \ldots \vee u_n$, where
all $u_i$ are equalities, are called {\it pseudo-equalities (or
pseudo-identities)}. The corresponding $\mathfrak X=T^*$ are
called pseudo-varieties. Action-type pseudo-varieties are given by
the formulas of the $w_1 \equiv 0 \vee \ldots \vee w_n \equiv 0 $
kind.

3. Formulas of the kind $u_1 \wedge u_2 \wedge \ldots \wedge u_n
\Longrightarrow u$, where all $u_i$ are equalities, are called
{\it quasi-equalities (or quasi-identities)}. Formulas $w_1 \equiv
0 \wedge \ldots \wedge w_n \equiv 0 \Longrightarrow w \equiv 0$ is
an action-type quasi-identity. We have quasi-varieties and
action-type quasi-varieties.

4. Universal formulas have the form $u_1 \vee  \ldots \vee u_n
\vee \neg v_1 \vee \ldots \vee \neg v_m$ where $u_i$ and $v_j$ are
equalities. Universal action-type formulas have the form $w_1
\equiv 0 \vee \ldots \vee w_n \equiv 0 \vee w_1 ' \not\equiv 0
\vee \ldots \vee w_m ' \not\equiv 0$. The corresponding classes
$\mathfrak X$ are universal classes and action-type universal
classes. We consider also classes $\mathfrak X$, axiomatizable by
arbitrary action-type formulas.

 A class $\frak X$ we call a
{\it saturated one}  if the inclusion $(V,G) \in \mathfrak X$
holds if and only if for the corresponding faithful representation
$(V, \overline G)$ holds $(V, \overline G)\in \mathfrak X$.

A class $\frak X$ we call a {\it right-hereditary one}, if $(V,G)
\in \mathfrak X \Longrightarrow (V,H) \in \mathfrak X$ holds for
every subgroup $H$ in $G$.

Finally, class of representations $\mathfrak X$ is called a {\it
right-local one}, if $(V,G) \in \mathfrak X$  holds if all $(V,H)$
belong to $\mathfrak X$, where  $H$ are finitary-generated
subgroups in $G$.

\section{Action-type axiomatizable classes}

\subsection{Theorem}

\begin{thm}
If class $\frak X$ is given by action-type formulas, then such
class is saturated, right-hereditary and right-local.
\end{thm}

Proof. Check first that the class is saturated.

Take a representation $(V,G)$ and the corresponding faithful
representation $(V, \overline G)$. Let $\beta _0 : G \to \overline
G$ be the natural homomorphism. Then for any $a \in V$  and $g \in
G$ we have $a \circ g = a \circ g ^{\beta_0}$. Consider, further,
homomorphisms

$$
\mu_0 = (\alpha, \beta) : W \to (V,G)
$$
and
$$
\mu = (\alpha, \beta \beta_0) : W \to (V,\overline G).
$$
We have the diagram

$$
\CD
W @>\mu_0>> (V,G)\\
@. @/SE/ \mu  // @VV \nu=(1,\beta_0) V\\
@.(V,\overline G)\\
\endCD
$$


Here $\mu = \mu_0 \nu$ and every $\mu$ can be represented in such
a way.

Note that the homomorphism $\beta: F \to G$ induces the
homomorphism $\beta : KF \to KG$ of group algebras. We have also
$\beta_0 : KG \to K \overline G$, and $a \circ h = a \circ
h^{\beta_0}$ for every $a \in V$ and $h \in KG$.

Take now an element $w = x_1 \circ u_1 + \ldots  + x_n \circ u_n$.
Then
$$
w^\alpha = x_1^\alpha \circ u_1^\beta + \ldots  + x_n^\alpha \circ
u_n^\beta= x_1^\alpha \circ u_1^{\beta\beta_0} + \ldots  +
x_n^\alpha \circ u_n^{\beta\beta_0}.
$$
From this follows that $w^\alpha\equiv 0$ holds in $(V,G)$ if and
only if the same holds in $(V,\overline G)$. In other words,
$$
\mu = (\alpha, \beta\beta_0) \in Val_{(V,\overline G)} (w\equiv 0)
\Leftrightarrow \mu_0 =(\alpha, \beta) \in Val_{(V,G)}(w\equiv 0)
$$
holds for every formula of the type $w\equiv 0$. Prove by
induction, that this is true for every action-type formula. Let it
be true for action-type formulas $u$ and $v$. Pass to $u \vee v$,
$u \wedge v$, $\neg u$, $\exists x \ u$. Let $\mu \in Val_{(V,
\overline G)} (u \vee v) = Val_{(V, \overline G)} (u) \vee
Val_{(V, \overline G)} (v)$ and let, say, $\mu \in Val_{(V,
\overline G)} (u)$. Then $\mu_0 \in Val_{(V, G)} (u)$, $\mu_0 \in
Val_{(V, G)} (u \vee v)$. Similarly we check that $\mu_0 \in
Val_{(V, G)} (u \vee v)\Longrightarrow \mu \in Val_{(V, \overline
G)} (u \vee v)$. Just the same for the case $u \wedge v$.

Let now $\mu \in  Val_{(V, \overline G)} (\neg u) = \neg
Val_{(V,\overline G)} (u)$, i.e., $\mu \not\in  Val_{(V, \overline
G)} ( u)$. Hence, $\mu_0 \not\in Val_{(V, G)} (u)$, $\mu_0 \in
Val_{(V, G)} (\neg u)$. The same for the opposite direction.

Consider now $\exists x u$. Let $\mu \in  Val_{(V, \overline G)}
(\exists x u) = \exists x  Val_{(V, \overline G)} ( u)$. Find $\nu
= (\alpha_1, \beta\beta_0) \in   Val_{(V, \overline G)} ( u)$,
$\alpha_1(x_1) = \alpha(x_1), x_1 \neq x$. Let us note here, that
the induction proceeds by any pairs $\mu$ and $\mu_0$ fitting the
diagram. Take $\nu_0 = (\alpha_1, \beta)$, $\nu = \nu_0
(1,\beta_0$). Inclusion $\nu \in Val_{(V, \overline G)} ( u)$ is
equivalent to the inclusion $\nu_0 \in Val_{(V,  G)} ( u)$. Since
$\alpha_1(x_1)=\alpha (x_1)$ for $x_1\neq x$, we conclude that
$\mu_0 \in  Val_{(V,  G)} (\exists x u)$. Similarly, $\mu_0 \in
Val_{(V, \overline G)}$ implies $\mu \in  Val_{(V, \overline G)}$.

Assume now that the class $\frak X$ is given by some set of
action-type formulas. Each of these formulas holds in the
representation $(V,G)$ if and only if it holds in the
representation $(V,\overline G)$. Hence, the class $\frak X$ is
saturated.

Check further that the class is right-hereditary. Take $(V,G)$ and
$(V,H)$, $H \subset G$. For every action-type formula $u $ verify
that
$$
Val_{(V,H)} (u) = Val_{(V,G)} (u) \cap Hom(W, (V,H)).
$$
Again we apply induction. The case when $u$ is an equality
$w\equiv 0$ is easily checked. It is also easily verified that if
$u$ and $v$ meet the condition, then the same is true for $u \vee
v$, $u \wedge v$ and $\neg u$. It is left to check $\exists x u$
if there is an equality for $u$. Take $\mu = (\alpha, \beta) \in
Val_{(V,H)} (\exists x u) = \exists x Val_{(V,H)}(u)$. Then $\mu
\in Hom(W, (V,H))$ and there is an element $\nu = (\alpha_1,
\beta)$ with $\alpha_1(x_1) = \alpha(x_1)$ if $x_1 \neq x$, $\nu
\in Val_{(V,H)}(u)$. By the condition, $\nu \in Val_{(V,G)}(u)$
and $\mu \in \exists x Val_{(V,G)} = Val_{(V,G)} \exists x u$,
$\mu \in Val_{(V,G)}\exists x u \cap Hom(W, (V,H))$. Conversely,
let $\mu \in Val_{(V,G)}(\exists x u) \cap Hom(W, (V,H))$, $\mu =
(\alpha, \beta)$ and $\beta$ is a homomorphism $F \to H$. Find
$\nu = (\alpha, \beta) \in Val_{(V,G)}( u) \cap Hom(W, (V,H)) =
Val_{(V,H)} (u)$, $\alpha_1(x_1)= \alpha(x_1), \ x_1 \neq x$. Now
$\mu \in Val_{(V,H)}(\exists x u) $. The equality is checked.
Assume now, that the action-type formula $u$ holds in the
representation $(V,G)$, $Val_{(V,G)}(u) = Hom(W, (V,G))$. We have
$Val_{(V,H)}( u) = Hom(W, (V,G)) \cap Hom(W, (V,H))= Hom(W,
(V,H))$, and $u$ holds in $(V,H)$. Thus, the class $\frak X$ is
right-hereditary.

Now it is only left to prove that the class $\frak X$ is
right-local. As earlier, assume that $\mathfrak X$ is given by
action-type formulas. We need to show, that if $u$ is action-type
formula, $(V,G)$ is a representation, $u$ holds on every $(V,H)$
where $H$ is finite-generated  subgroup in $G$, then $u$ holds in
$(V,G)$ as well.

We need some auxiliary material. Let a homomorphism $\mu =
(\alpha, \beta): W \to (V,G)$, $W = W(X,Y) = (XKF,F)$ be given.
Take a subset $Y_0$ in $Y$ and change $\beta$ by $\beta ' : F \to
G$, coinciding with $\beta$ on $Y_0$ and sending all the rest
$y$'s to the unit. Take further $\mu ' = (\alpha ', \beta ')$,
where $\alpha '$ coincides with $\alpha$ on $\mathfrak X$. Take
now an arbitrary $w = x_1 \circ u_1 + \ldots x_n \circ u_n$ such
that supports of all $u_i$ belong to $Y_0$ (in other words all
elements $u_i$ are expressed via variables from $Y_0$). Then

$$
w^{\alpha '} = x_1^{\alpha '} \circ u_1^{\beta '} + \ldots
x_n^{\alpha '} \circ u_n^{\beta '} = x_1^{\alpha } \circ
u_1^{\beta} + \ldots x_n^{\alpha} \circ u_n^{\beta} = w^{\alpha}.
$$

For an arbitrary action-type formula $u$ consider its $Y$-support
$\Delta_Y (u)$. We have:

$$
\Delta_Y(u) = \Delta_Y (\neg u) = \Delta_Y (\exists x u),
$$

$$
\Delta_Y(u \vee v) = \Delta_Y ( u \wedge v) = \Delta_Y( u)\cup
\Delta_Y( v).
$$

Let now $\mu = (\alpha,\beta)$ be a homomorphism and $u$ an
action-type formula, $Y_0 \supset \Delta_Y(u)$, $Y_0$ be a finite
set. Pass to $\mu ' = (\alpha ', \beta ')$ by $Y_0$.  The
following property always takes place in these conditions:

$$
\mu = (\alpha, \beta) \in Val_{(V,G)}(u) \Leftrightarrow \mu ' =
(\alpha ', \beta ') \in Val_{(V,G)}(u).
$$

Let us check it. The property is evident for identities. Now let
it be true for some action-type formula $u$. Pass to $\neg u$ and
$\exists x u$.  We have $\Delta_Y (\neg u) = \Delta_Y (\exists x
u) = \Delta_Y(u)$. Take an arbitrary finite $Y_0$, including these
supports. Let

$$\mu  \in Val_{(V,G)}(\neg u), \ i.e., \  \mu \not
\in Val_{(V,G)}(u).$$

Then $\mu ' \not \in Val_{(V,G)}(\neg u)$ and $\mu ' \in
Val_{(V,G)}(\neg u)$. The opposite is checked similarly.

Pass to $\exists x u$ with the same $Y_0$. Given $\mu = (\alpha,
\beta) \in Val_{(V,G)}(\exists x u) = \exists x Val_{(V,G)}(u)$,
we have $\mu ' = (\alpha ', \beta ')$ by $Y_0$. Besides that, we
can select $\mu_1 = (\alpha_1, \beta)$ such that $\mu_1 = \in
Val_{(V,G)}(u)$, $\alpha_1 (x_1) = \alpha (x_1)$ for $x_1 \neq x$.
Take once more $\mu_1 ' = (\alpha_1 ', \beta ')$ for $\mu_1 =
(\alpha_1, \beta)$ by $Y_0$. Now $ \mu_1 \in Val_{(V,G)}(u)$
implies $ \mu_1 '  \in Val_{(V,G)}(u)$. Recall that $\alpha$ and
$\alpha '$ coincide on the set $X$, as well as $\alpha_1$ and
$\alpha_1 '$. Then for $x_1 \neq x$ we have $\alpha_1 ' (x_1) =
\alpha_1 (x_1) = \alpha (x_1)\alpha ' (x_1)$. Comparing $\mu '$
with $\mu_1 '$ we conclude: $\mu ' \in Val_{(V,G)}(\exists x u)$.
Similarly we derive $\mu ' \in Val_{(V,G)}(\exists x u)
\Rightarrow \mu  \in Val_{(V,G)}(\exists x u)$.

Let now the property holds for $u$ and $v$. Check $u \vee v$, $u
\wedge v$. Take $Y_0 \supset \Delta_Y(u \vee v) = \Delta_Y(u
\wedge v) = \Delta_Y(u)\cup \Delta_Y(v)$. Take $\mu = (\alpha,
\beta)$ and pass to $\mu ' = (\alpha ', \beta ')$ by $Y_0$.

Given $\mu  \in Val_{(V,G)}(u \vee v) = Val_{(V,G)}(u) \cup
Val_{(V,G)}(v)$. Let $\mu  \in Val_{(V,G)}(u)$. Since $Y_0 \supset
\Delta_Y(u)$, then $\mu ' \in Val_{(V,G)}(u)$ and $\mu ' \in
Val_{(V,G)}(u \vee v)$.The rest is done in the same way. The
property is checked for every $u$.

We continue the proof of the theorem. Let the formula $u$ be hold
in every $(V,H)$. We have: $Val_{(V,H)}(u) = Val_{(V,G)}(u) \cap
Hom(W, (V,H))$. Take an arbitrary $\mu = (\alpha, \beta)$ and pass
to $\mu '= (\alpha ', \beta ')$ by $u$. Take $Im \ \beta '$ as
$H$. Then $\mu ' \in Hom(W, (V,H)) = Val_{(V,H)}(u)$. But then
$\mu ' \in Val_{(V,G)}(u)$ and this gives $\mu \in
Val_{(V,G)}(u)$. This is true for every $\mu$ and $Val_{(V,G)}(u)
= Hom(W, (V,G))$. The formula $u$ holds in the representation
$(V,G)$.

Let now the class $\frak X$ is given by the set $T$ of action-type
formulas and for the representation $(V,G)$ every representation
$(V,H)$ belongs to $T*$, i.e., satisfies all $u \in T$. Then all
these $u$ are fulfilled in $(V,G)$ and $(V,G)$ belongs to
$\mathfrak X$. Hence, the goal local property takes place and the
theorem is proved.

There arises the following main

\begin{prob}

Is it true that an axiomatized class $\frak X$ is action-type
axiomatized if and only if it is saturated, right-hereditary and
right-local?

\end{prob}

Theorem 2.1 solves the question in one direction, but the opposite
claim  requires additional considerations.

\subsection{Theorem}

\begin{thm}
 For varieties,
pseudovarieties, quasivarieties and universal classes of
reperesentations  Problem 1 has a positive solution.
\end{thm}

Proof.

The case of varieties is studied in \cite{PV}. Apply the general
approach, also outlined in \cite{PV}.

Let $\mathfrak X$ be a saturated class and $G$ a group. There are
representations with the acting group $G$ in $\mathfrak X$. Denote
by $\mathfrak X_G$ a {\it $G$-layer} in $\mathfrak X$, all $(V,G)$
in $\mathfrak X$ are with the given group $G$.

Let now $\mathfrak X$ be a pseudovariety satisfying the closure
conditions from Problem 1. We need to show that $\mathfrak X$ is
given by action-type pseudoidentities. Show first that the class
$\frak X$ with these closure conditions is a pseudovariety if and
only if every layer $\mathfrak X_G$ is a pseudovariety of
$G$-modules.

Recall that the class of algebras $\Theta$ is a pseudovariety if
and only if it is closed in respect to ultra-products, subalgebras
and homomorphic images (for one-sorted algebras see \cite{Ma}, the
generalization for multi-sorted ones see in \cite{Gva}).

Let now $\mathfrak X$ be a pseudovariety of representations and
let this class be is saturated, right-hereditary and right-local .
Pass to the layer $\mathfrak X_G$. It is evident that it is closed
under subalgebras and homomorphisms and we need to check that it
is closed under ultra-products.

Let $I$ be a set and $D$ a filter of subsets in $I$. We have a
representation $(V_\alpha, G)$ for every $\alpha \in I$. Take a
cartesian product $((\prod_\alpha V_\alpha), G^I)$ and let
$(V_0,H)$ be a congruence determined by the filter $D$. The
corresponding filtered product is $((\prod_\alpha V_\alpha )/ V_0,
G^I / H)$. It belongs to the class $\frak X$ if $D$ is an
ultrafilter. Since $\mathfrak X$ is saturated, the representation
$((\prod_\alpha V_\alpha )/ V_0, G^I)$ also belongs to $\mathfrak
X$. Embed now $G \to G^I$ as a diagonal. Since $\mathfrak X$ is
hereditary, the $((\prod_\alpha V_\alpha) / V_0, G)$ belongs to
$\mathfrak X$ and, consequently, to $\mathfrak X_G$. It is easy to
understand, that this representation is isomorphic to
ultra-product by the filter $D$ of all $(V_\alpha, G)$ with the
fixed $G$. Hence,  all $\mathfrak X_G$ are pseudovarieties.

Let now the class $\frak X$ be saturated, right-hereditary and all
$\mathfrak X_G$ be pseudovarieties. Show that $\mathfrak X$ are
pseudovarieties. Check first, that the class $\frak X$ is
hereditary and closed under homomorphisms. Take $(V,G)\in
\mathfrak X$ and let $(V_0,H)$ be a subrepresentation. From
right-hereditaryty follows that $(V,H)$ belongs to $\mathfrak X$.
Consider further a  layer $\mathfrak X_H$. It is a pseudovariety
and, consequently, $(V_0,H)\in \mathfrak X_H$ and $(V_0,H)\in
\mathfrak X$.

Let now a surjective homomorphism $\mu = (\alpha, \beta): (V,G)\to
(V_1,G_1) $, $(V,G)\in \mathfrak X$ be given. The homomorphism
$\beta : G \to G_1$ determines the representation $(V_1,G)$. We
have a homomorphism $(\alpha, 1): (V,G) \to (V_1,G)$. Since the
layer $\mathfrak X_G$ is closed under homomorphisms, then $(V_1,G)
\in \mathfrak X_G$ and $(V_1,G)\in \mathfrak X$. The class $\frak
X$ is saturated. This means that $(V_1,G_1) \in \mathfrak X$.

It is left to check that the class $\frak X$ is ultra-closed. Let
$I$ be a set, $D$ ultrafilter on $I$. Given $(V_\alpha, G_\alpha)
\in \mathfrak X$, $\alpha \in I$, take $G = \prod_\alpha
G_\alpha$. We have $\pi_\alpha : G \to G_\alpha$ determining the
representation $(V_\alpha, G)$. It belongs to $\mathfrak X$ since
$\mathfrak X$ is saturated. Consider the layer $\mathfrak X_G$ and
take in it an ultra-product by $D$. We have: $((\prod V_\alpha) /
D, G )\in \mathfrak X_G$ and hence it belongs to $\mathfrak X$.
Ultra-product by $D$ of all $(V_\alpha, G_\alpha)$ is a
homomorphic image of the representation $((\prod V_\alpha) /
D,G)$. Since we have checked closure in respect to homomorphisms,
then the corresponding ultra-product contains in $\mathfrak X$.
Hence, $\mathfrak X$ is a pseudovariety with the given closure
conditions.

Let further $\mathfrak X$ be a pseudovariety. Prove that
$\mathfrak X$ is determined by action-type pseudo-identities. Take
a countable set $Y$ and a free group $F=F(Y)$. Consider a layer
$\mathfrak X_F$. It is a pseudovariety and it is determined by
pseudo-identities. Take the corresponding set of pseudo-identities
and let it be indexed by the set $I$. We identify $I$ with that
set.

Consider pseudo-identities $$w_1^k \equiv 0 \vee \ldots w_{n_k}
^k\equiv 0,\ \ k \in I.$$ Denote by $u^k$ the given $k$-th
pseudo-identity, $u^k = u^k(x_1, \ldots, x_n; y_1, \ldots, y_m)$.
Besides, we have $w_i^k = w_i^k(x_1, \ldots, x_n; y_1, \ldots,
y_m)$. We consider all $x_1, \ldots, x_n$ as variables while $y_1,
\ldots, y_m$ are constants. Show that the same set $I$ with the
varying $y$'s determines the class $\frak X$. Take a
representation $(V,G)$ and show that it belongs to the class
$\frak X$ if and only if all pseudo-identities of the given set
hold in it.

Let $(V,G) \in \mathfrak X$. Consider a homomorphism $\mu =
(\alpha, \beta): W = W(X,Y) = (XKF,F) \to (V,G)$. Using $\beta: F
\to G$, determine a representation $(V,F)$. Since the class $\frak
X$ is saturated and right-hereditary, we conclude that $(V,F) \in
\mathfrak X$. Then $(V,F) \in \mathfrak X_F ,$ and the set $I$
holds in $(V,F)$. Hence, for every $k \in I$ there is $w_i^k
=w_i^k(x_1, \ldots, x_n; y_1, \ldots, y_m)$ with
$w_i^k(x_1^\alpha, \ldots, x_n^\alpha; y_1, \ldots, y_m) = 0$. But
$w_i^k(x_1^\alpha, \ldots, x_n^\alpha; y_1, \ldots, y_m) =
w_i^k(x_1^\alpha, \ldots, n^\alpha; y_1^\beta, \ldots,
y_m^\beta)$. Thus, $u^k$ holds in $(V,G)$ under the homomorphism
$\mu = (\alpha, \beta)$. It is true for every $k$ and $\mu$. This
means that the set of formulas $I$ holds in the representation
$(V,G)$.

Conversely, let the set $I$ holds in $(V,G)$. We need to check
that $(V,G) \in \mathfrak X$. Take a finitely generated subgroup
$H$ in $G$ and a surjection $\beta : F \to H$. Take $(V,F)$
defined by this surjection. For an arbitrary $\alpha : \mathfrak X
\to V$ we have a $KF$-homomorphism $\alpha : XKF \to V$. We have
also $\mu = (\alpha, \beta): W \to (V,G)$. For every $k$ find
$w_i^k(x_1, \ldots, x_n; y_1, \ldots, y_m)$ with
$w_i^k(x_1^\alpha, \ldots, n^\alpha; y_1^\beta, \ldots, y_m^\beta)
= 0 = w_i^k(x_1^\alpha, \ldots, n^\alpha; y_1, \ldots, y_m)$. This
means that the set $I$ with the fixed $y$'s  is hold in $(V,F)$,
i.e., $(V,F)\in \mathfrak X_F$. But then $(V,F) \in \mathfrak X$
and $(V,H) \in \mathfrak X$. This holds for every finitely
generated subgroup $H$ in $G$ and, therefore, $(V,G)\in \mathfrak
X$.

The case of pseudovariety is  completed. The reasoning of the same
type is used for quasi-varieties and universal classes.

Consider quasivarieties in more detail, taking into account their
special role in algebraic geometry \cite{PTs}.

It is known (see \cite{Ma}, \cite{GL}) that if $\mathfrak X$ is a
class of algebras and $qVar(\mathfrak X)$ is a quasivariety
generated by the class $\frak X$, then

$$qVar(\mathfrak X) = SCC_{up} (\mathfrak X).$$

Here $S$ is an operator of taking of subalgebras, $C$ cartesian
products and $C_{up}$ an operator of taking  ultra-products. This
implies also that it is enough to check that the corresponding
class is closed in respect to the operator $S$ and arbitrary
filtered products. Besides, we need to assume that there are
one-element subalgebras in the class.

Using all this, prove the theorem for quasivarieties. Let us prove
that if $\mathfrak X$ is a saturated, right-hereditary and
right-local quasivariety, then this $\mathfrak X$ is an
action-type quasivariety.

As earlier, select the following two items of the proof. Check
first of all that if $\mathfrak X$ satisfies the pointed closure
conditions, then $\mathfrak X$ is a quasivariety if and only if
every layer $\mathfrak X_G$ is a quasivariety. Then using the
layer $\mathfrak X_F$ and closure conditions, prove that the class
$\mathfrak X$ is determined by action-type quasivarieties. The
first item of the proof copies the one for pseudo-varieties with
the only difference that we work with filtered products instead of
cartesian and ultra-products. The second item follows the idea
used for pseudo-varieties. The same reasoning with two items we
use for universal classes. Now the theorem 2.2 is proved.

\section{General theorem}

\subsection{Some general considerations}

In the previous theorem we used the fact that if the class $\frak
X$ is given by axioms of a special kind, then every layer
$\mathfrak X_G$ is determined by the axioms of the same kind as
well. However we cannot say anything similar if the type of the
axioms is not fixed. We cannot claim, that if the class $\frak X$
is axiomatized and meets the demanded closure conditions, then the
layer $\mathfrak X_F$ is also axiomatized. We will study this
obstacle later on.

Let us note that it is convenient to use the following
terminology. We call a class $\frak X$ {\it strictly saturated} if
it is saturated, right-hereditary and right-local.

\subsection{The theorem}

\begin{thm}
Let $\mathfrak X$ be an axiomatized class of representations and a
layer $\mathfrak X_F$ also axiomatized. Such $\mathfrak X$ is
action-type axiomatized if and only if it is strictly saturated.
\end{thm}

Proof. Let us make some additional remarks on the layer $\mathfrak
X_F$. First of all, we denote by $Rep-KF$ the class of all
$K$-representations of the free group $F$. The class $Rep-KF$ of
all $(V,F)$ is a subvariety in $Rep-K$.

Every representation $(V,F)$ we consider now as a one-sorted
algebra, whose morphisms have the type $\alpha= (\alpha,1): (V,F)
\to (V_1, F)$, and they commute  with the action of the group $F$.
The free representation in $Rep-KF$ is the same $(XKF,F)$. Logic
of the representation $Rep-KF$ consists of action-type formulas
but we consider these formulas from another perspective: we view
variables of the set $Y$ as constants, they are immutable. If $u$
is an action-type formula in the logic over $Rep-K$, then the same
$u$ in the logic $Rep-KF$ is denoted by $u^0$. If $T$ is the set
of action-type formulas, then correspondingly we consider $T^0$.
Let the set $T$ determines a class $\frak X$. Whether we can claim
that $T^0$ determines the layer $\mathfrak X_F$? Probably not. We
can claim that every $u^0 \in T^0$ holds in every $(V,F) \in
\mathfrak X_F$, but we cannot claim that if $T^0$ holds in some
$(V,F)$, then $T$ holds in $(V,F) \in \mathfrak X$ as well. There
are no general arguments for the feature of the layer $\mathfrak
X_F$ to be axiomatizable.

Let us pass to the proof of the theorem 3.2. Let the conditions of
this theorem hold for $\mathfrak X$ and the layer $\mathfrak X_F$
is axiomatized by a set of formulas $T^0$ where $T$ is a set of
action-type formulas. Show that this $T$ determines the class
$\frak X$. We use the same reasoning as before.

Take a representation $(V,G)$. Let it belong to $\mathfrak X$.
Check that all formulas from $T$ hold in $(V,G)$.

Consider a homomorphism $$\mu = (\alpha, \beta): (XKF,F) \to
(V,G).$$

Applying $\beta: F \to G$, pass to the representation $(V,F)$. Let
$Im \ \beta = H$. We have $(V,H)$, and this representation belongs
to $\mathfrak X$ since it is right-hereditary. The class $\frak X$
is saturated, hence the representation $(V,F)$ belongs to
$\mathfrak X$. But then it belongs to the layer $\mathfrak X_F$
and all axioms of the set $T^0$ hold in it.

Here we deviate a bit from the main stream and consider the
following property.

Given a homomorphism $\mu = (\alpha, \beta): (XKF,F) \to (V,G)$,
pass, as earlier, to $(V,F)$. Then $\alpha \in Val_{(V,F)} (u^0)
\Leftrightarrow \mu \in Val_{(V,G)} (u) $ for every action-type
formula $u$. Prove this by induction. Let first $u$ be $w \equiv
0$, $w = x_1 \circ u_1 + \ldots + x_n \circ u_n$. Then $w^\alpha =
x_1^\alpha \circ u_1 + \ldots + x_n^\alpha \circ u_n = x_1^\alpha
\circ u_1^\beta + \ldots + x_n^\alpha \circ u_n^\beta$. Thus,
$w^\alpha = 0$ holds in $(V,F)$ if and only if the same holds in
$(V,G)$. The property is checked for elementary formulas.

Let the property be checked for some $u$. Verify it for $\neg u$
and $\exists x u$. Take $\alpha \in Val_{(V,F)} (\neg u^0)$, i.e.,
$\alpha \not \in Val_{(V,F)} (u^0)$, and $\mu \not \in Val_{(V,G)}
(u)$, i.e., $\mu \in Val_{(V,G)} (\neg u)$. The same for the
opposite direction.

Let now $\alpha \in Val_{(V,F)} (\exists x u^0) = \exists x
Val_{(V,F)} (u^0)$. Find $\alpha ' \in Val_{(V,F)} (u^0)$ with
$\alpha '(x_1) =\alpha (x_1)$ for $x_1 \neq x$. We have $\mu =
(\alpha, \beta)$. Take $\mu ' = (\alpha ', \beta) = (\alpha ',1)
(1, \beta)$. Since $\alpha ' \in Val_{(V,F)} (u^0)$ and $(V,F)$ is
built from $(V,G)$ by $\beta$, then $\mu ' \in Val_{(V,G)} (u)$.
Together with $\alpha '(x_1)=\alpha (x_1)$ this gives $\mu \in
Val_{(V,G)} (\exists x u)$. Let $\mu \in Val_{(V,G)} (\exists x
u)$, $\mu = (\alpha, \beta)$. Take $\mu ' = (\alpha ', \beta)\in
Val_{(V,G)} (u)$, $\alpha '(x_1) =\alpha (x_1)$, $x_1 \neq x$. We
have $\alpha ' \in Val_{(V,F)} (u^0)$ and $\alpha \in Val_{(V,F)}
(\exists x u^0)$. The cases $u \vee v$ and $u \wedge v$ are
checked in the same way. The property is checked.

Let us return to the proof of the theorem. Let now $u \in T$, $u^0
\in T^0$. Then the formula $u^0$ holds in the representation
$(V,F)$, and for every $\alpha : (XKF,F) \to (V,F)$ we have
$\alpha \in Val_{(V,F)} (u^0)$. Apply this to the initial $\mu =
(\alpha, \beta)$. Then $\mu \in Val_{(V,G)} (u)$. This holds for
every $\mu$ and $Val_{(V,G)} (u) = Hom(W,(V,G))$. The formula $u$
holds in $(V,G)$.

Let now the set $T$ holds in $(V,G)$. We need to check that $(V,G)
\in \mathfrak X$. Since $\mathfrak X$ is right-local, it is enough
to verify that $(V,H) \in \mathfrak X$ for every finitely
generated subgroup $H \subset G$.

Consider such $H$ and a surjection $\beta F \to H$. Take $(V,F)$
by it. For an arbitrary $\alpha X \to V$ and pass to
$KF$-homomorphism $\alpha : XKF \to V$. We have also $\mu =
(\alpha, \beta): W \to (V,G)$.

As we have seen, $\alpha \in Val_{(V,F)} (u^0) \Leftrightarrow \mu
= (\alpha, \beta) \in Val_{(V,G)} (u)$ holds for every formula
$u^0 \in T^0$, $u \in T $. Since $\mu \in Val_{(V,G)} (u)$ for
every $\mu$, then $\alpha \in Val_{(V,F)} (u^0)$ for every $\alpha
: XKF \to V$. This means that the set $T^0$ holds in $(V,F)$.
Since $T^0$ determines the layer $\mathfrak X_F$, then $(V,F) \in
X_F$, and, hence, $(V,F) \in \mathfrak X$. This implies that
$(V,H) \in \mathfrak X$ and, therefore, we have $(V,G) \in
\mathfrak X$.

The theorem is proved in one direction, and what's more, the
initial class $\frak X$ need not be axiomatizable. The opposite
direction follows from the theorem 2.1, and here we do not demand
the layer $\mathfrak X_F$ to be axiomatizable.

The main Problem 1 remains open in full generalality.



\end{document}